\documentstyle[11pt]{article}

\makeatletter

\@addtoreset{equation}{section}
\makeatother
\def\supp{{\rm supp}\ }

\def\<{\langle}
\def\>{\rangle}

\newtheorem{lem}{Lemma}[section]
\newtheorem{theo}{Theorem}[section]
\newtheorem{rem}{Remark}[section]

\makeatletter
   
   \@addtoreset{equation}{section}
\makeatother

\begin{document}
\title{\bf Fast Energy Decay for Wave Equations with Variable Damping Coefficients in the $1$-D Half Line}
\author{Ryo IKEHATA\thanks{Corresponding author: ikehatar@hiroshima-u.ac.jp} and Takeshi KOMATSU\\ {\small Department of Mathematics, Graduate School of Education, Hiroshima University} \\ {\small Higashi-Hiroshima 739-8524, Japan}}
\maketitle
\begin{abstract}
We derive fast decay estimates of the total energy for wave equations with localized variable damping coefficients, which are dealt with in the one dimensional half line $(0,\infty)$. The variable damping coefficient vanishes near the boundary $x = 0$, and is effective critically near spatial infinity $x = \infty$. 
\end{abstract}

\section{Introduction}
\footnote[0]{Keywords and Phrases: Localized damping; Wave equation; Multiplier method; Total energy decay; Weighted initial data.}
\footnote[0]{2010 Mathematics Subject Classification. Primary 35L20; Secondary 35L05, 35B33, 35B40.}
We consider the $1$-dimensional initial-boundary value problem for linear dissipative wave equation in the half line $(0,+\infty)$:
\begin{equation}
u_{tt}(t,x) - u_{xx}(t,x) + V(x)u_{t}(t,x) = 0,\ \ \ (t,x)\in (0,\infty)\times (0,\infty),
\end{equation}
\begin{equation}
u(0,x)= u_{0}(x),\ \ u_{t}(0,x)= u_{1}(x),\ \ \ x\in (0,\infty),
\end{equation}
\begin{equation}
u(t,0)=0,\ \ t\in (0,\infty),
\end{equation}
where the initial data $(u_{0},u_{1})$ are compactly supported with the size $R > 0$:
$$u_{0} \in H^{1}_{0}(0,\infty),\ \ u_{1} \in L^{2}(0,\infty),\ \ \supp u_{i}\subset [0,R], (i=0,1).$$

Let us first talk about several related results of the Cauchy problem in ${\bf R}^{N}$ of the equation:
\begin{equation}
u_{tt}(t,x) - \Delta u(t,x) + V(x)u_{t}(t,x) = 0.
\end{equation}
In the case when $V(x) = $ constant $> 0$, Matsumura \cite{Ma} derived a historical $L^{p}$-$L^{q}$ decay estimates of solutions of (1.4). Later on more detail investigations such as the asymptotic profiles of solutions were studied by Nishihara \cite{Ni} and the references therein.

In the case of mixed problems of (1.4) in the exterior domain $\Omega \subset {\bf R}^{N}$ Nakao \cite{N} and Ikehata \cite{ike-2} have derived precise decay estimates of the total energy for the initial data belonging to the energy class and weighted $L^{2}$-class, respectively. In these cases they have dealt with the localized damping coefficient $V(x)$, which is effective only near infinity such as $V(x) > V_{0}$ for $\vert x\vert \gg 1$, and the star-shaped boundary $\partial\Omega$. In particular, fast decay estimates of the total energy like $E(t) = O(t^{-2})$ was obtained in \cite{ike-2} by a special multiplier method developed in \cite{IM}, which modified the Morawetz one \cite{Mora}. Racke \cite{R} also derived the $L^{p}$-$L^{q}$ estimates for solutions to (1.4) with constant $V(x)$ in the exterior domain case by relying on so called the generalized Fourier transform (in fact, more general form of equations and systems are considered). For related results in the case when the damping is effective near infinity one can cite the papers written by Aloui-Ibrahim-Khenissi \cite{AIK} and Daoulatli \cite{D}, where they studied the fast decay of the total and local energy under the so called GCC assumption due to Lebeau \cite{L}. Exponential decay of the total energy for wave equations with mass (i.e., Klein-Gordon type) and localized damping terms was investigated by Zuazua \cite{Z} in the exterior domain case. Furthermore, by Watanabe \cite{W} an effective application of the method developed in \cite{ike-2} is recently introduced, and there a global existence and decay estimates of solutions to a nonlinear Cauchy problem of the quasilinear wave equation with a localized damping near infinity are studied. 

On the other hand, recently in \cite{ike-4} the critically decaying damping coefficient $V(x)$ near spatial infinity was considered for (1.4) in the whole space ${\bf R}^{N}$, in fact, they have studied so called the critical case with $\delta = 1$ below:
\begin{equation}
\frac{V_{0}}{(1+\vert x\vert)^{\delta}}\leq V(x) \leq \frac{V_{1}}{(1+\vert x\vert)^{\delta}}\quad x \in {\bf R}^{N},
\end{equation}
and they have announced the fact that $E(t) = O(t ^{-N})$ if $V_{0} > N$, while if $V_{0} \leq N$, then (roughly speaking) $E(t) = O(t ^{-V_{0}})$.  In this connection, we call it sub-critical damping in the case of (1.5) with $\delta \in [0,1)$, while in the case when $\delta > 1$, (1.5) is called as super-critical damping.  The precise decay order of several norms of solutions for the Cauchy problem in ${\bf R}^{N}$ of (1.4) with (1.5) were investigated by Todorova-Yordanov \cite{TY} in the sub-critical case. Furthermore, in the super critical case it is well-known by Mochizuki \cite{Mochi} that the solution of the Cauchy problem of (1.4) is asymptotically free.

While, in the constant damping case with $V(x) = $ constant of (1.1)-(1.3) (i.e., problem in the half line) was studied in \cite{ike-1}, and he derived $E(t) = O(t ^{-2})$ ($t \to +\infty$).  In the localized damping case, it seems that there are no any previous research papers that are dealing with the half space problem and critical damping near infinity. In this connection, one notes that in \cite{ike-3} Ikehata-Inoue have derived $E(t) = O(t^{-1})$ as $t \to +\infty$ for $V_{0} > 1$ (see also Mochizuki-Nakazawa \cite{MN} and Uesaka \cite{U} for the topics on the energy decay property). In fact, they treated a more general nonlinear equation like
\[u_{tt}(t,x) - \Delta u(t,x) + V(x)u_{t}(t,x) + \vert u(t,x)\vert^{p-1}u(t,x) = 0,\]
in the critical damping case of (1.5).  

So, a natural question arises whether one can derive fast decay result like $E(t) = O(t^{-2})$ in the case when $V(x)$ vanishes near the boundary and is effective critically (i.e., $\delta = 1$ in (1.5)) only near spatial infinity $x = \infty$. The main purpose of this paper is to answer this natural question by modifying a technical method originally developed in \cite{ike-3}.\\

Our assumptions on $V(x)$ are as follows:\\

({\bf A}) $V \in L^{\infty}(0,\infty) \cap C([0,\infty))$, and there exist three constants 
\[V_{0} > {\bf 2}, \quad V_{1} \in [V_{0},+\infty), \quad L_{2} > 0\]
such that
\[\frac{V_{0}}{1+x}\leq V(x) \leq \frac{V_{1}}{1+x},\quad x \in [L_{2}, +\infty),\]
\[0 \leq V(x) \quad (x \in [0,\infty) ).\]

Our new results read as follows.

\begin{theo}\,Let $V(x)$ satisfy the assumption $({\bf A})$. If the initial data $(u_{0},u_{1})\in H^{1}_{0}(0,\infty) \times L^{2}(0,\infty)$ further satisfies
$$supp \,u_{0} \cup supp\,u_{1} \subset [0,R]$$
for some $R > 0$, then for the solution $u \in C([0,\infty);H^{1}_{0}(0,\infty))\cap C^{1}([0,\infty);L^{2}(0,\infty))$ to problem $(1.1)-(1.3)$ one has
$$E(t) = O(t^{-2}),\ \ \ \ (t \rightarrow \infty).$$
\end{theo}
\begin{rem}{\rm In the $1$-D whole space case, if we consider the corresponding Cauchy problem, we can have
$$E(t) = O(t^{-1}),\ \ \ \ (t \rightarrow \infty)$$
provided that $V_{0} >1$ (cf., \cite{ike-4}, \cite{U}). So, the results seem to reflect the half space property itself. This discovery of the new number $V_{0} > {\bf 2}$  in the half space case is one of our maximal contribution. Even if one takes $L_{2} = 0$ formally, the obtained result in Theorem {\rm 1.1} is new. }

\end{rem}
\begin{rem}{\rm In the case when $V_{0} \leq {\bf 2}$, we still do not know the exact rate of decay for $E(t)$. However, from \cite[Theorem 4.1] {ike-4} one can state a conjecture that $E(t) = O(t^{-V_{0}})$ ($t \to +\infty$) if $V_{0} \leq 2$. Furthermore, the result will be generalized to the $N$-dimensional half space problems of (1.1)-(1.3). We can present a conjecture that if $V_{0} > N+1$, then $E(t) = O(t^{-(N+1)})$ ($t \to +\infty$) in the half space case.}
\end{rem}
\begin{rem}{\rm The assumption ({\bf A}) implies that the damping coefficient $V(x)$ can vanish near the boundary $x = 0$, so the damping is effective only near infinity. But, the effectiveness of the damping $V(x)$ near infinity corresponds to so called the critical case of (1.5) with $\delta = 1$ (see \cite{ike-4}). In this sense, the result obtained in Theorem 1.1 is derived under two types of difficult situations, that is, one is vanishing damping, and the other is critical one.}
\end{rem}

This paper is organized as follows. In section 2 we shall prove Theorem 1.1 by relying on a multiplier method which was originally introduced by the first author's collaborative work \cite{ike-3} and \cite{IM}. \\

{\bf Notation.}\,{\small Throughout this paper, $\| \cdot\|$ means the usual $L^2(\Omega)$-norm. The total energy $E(t)$ to the solution $u(t,x)$ of (1.1) is defined by
$$E(t) := \frac{1}{2}(\| u_t(t,\cdot)\|^2+\| u_{x}(t,\cdot)\|^2). 
$$
Furthermore, one sets as a $L^{2}(0,\infty)$ inner product:
\[(u,v) := \int_{0}^{+\infty}u(x)v(x)dx.\]
A weighted function space is defined as follows: $u \in L^{1,1/2}(0,\infty)$ iff $u \in L^{1}(0,\infty)$ and
$$\|  u \|_{1,1/2} := \int^{\infty}_{0}\sqrt{1+x}|u(x)|dx < +\infty .$$
Furthermore, one sets
\[X_{1}(0,\infty) := C([0,\infty);H^{1}_{0}(0,\infty))\cap C^{1}([0,\infty);L^{2}(0,\infty)).\]}

\section{Proof of Results}

We first prepare an important lemma, which is borrowed from \cite[Lemma 2.1]{ike-1}.
\begin{lem} It holds that
$$\sup_{0 \leq x < +\infty} (\frac{|u(x)|}{\sqrt{1+x}}) \leq \| u_{x} \|$$
for all $u \in H^{1}_{0}(0,\infty)$. 
\end{lem} 


\begin{lem} 
Let $u \in X_{1}(0,\infty)$ be the solution to problem $(1.1)-(1.3)$. Then it is true that
$$\frac{d}{dt}{\cal E}(t)+{\cal F}(t)=0,\quad (t \geq 0)$$
where 
$${\cal E}(t)=\frac{1}{2}\frac{d}{dt}\int^{\infty}_{0}[f(u_{t}^{2}+u_{x}^{2}) + 2guu_{t} + (gV - g_{t})u^{2} + 2hu_{x}u_{t}]dx,$$
$${\cal F}(t) = \frac{1}{2}\int^{\infty}_{0}(2fV - f_{t} - 2g + h_{x})u_{t}^{2}dx + \frac{1}{2}\int^{\infty}_{0}(2g - f_{t} + h_{x})u_{x}^{2}dx$$
$$ + \frac{1}{2}\int^{\infty}_{0}(g_{tt}-g_{t}V)u^{2}dx + \int^{\infty}_{0}(hV - h_{t})u_{x}u_{t}dx - \frac{1}{2}\int^{\infty}_{0}\frac{d}{dx}(hu_{x})dx.$$
Moreover, $f=f(t),\ \  g=g(t)\ \  and\ \  h=h(t,x)$ are all smooth functions specified later on.

\end{lem}
{\it Proof.}\, Since one considers the weak solutions $u(t,x)$ of the problem (1.1)-(1.3), we can assume to check the desired identity that the solution $u(t,x)$ is sufficiently smooth, and vanishes for large $x \gg 1$. \\
{\bf Step 1.} Multiplying the both sides of  (1.1) by $fu_{t}$ and rearranging it one can get:\\
$$fu_{t}u_{tt} - fu_{t}u_{xx} + fVu_{t}^{2}$$
$$= f\frac{1}{2}\frac{d}{dt}u_{t}^{2} - \frac{d}{dx}(fu_{t}u_{x}) + \frac{1}{2}\frac{d}{dt}fu_{x}^{2} - \frac{1}{2}f_{t}u_{x}^{2} + fVu_{t}^{2}$$
$$= \frac{1}{2}\frac{d}{dt}[f(u_{t}^{2} + u_{x}^{2})] + (fV - \frac{1}{2}f_{t})u_{t}^{2} - \frac{1}{2}f_{t}u_{x}^{2} - \frac{d}{dx}(fu_{t}u_{x}) = 0.$$
{\bf Step 2.} Multiplying the both sides of  (1.1) by $gu$ and rearranging it one has:\\
$$guu_{tt} - guu_{xx} + gVu_{t}u$$
$$= g\frac{d}{dt}uu_{t} - gu_{t}^{2} - \frac{d}{dx}(guu_{x}) + gu_{x}^{2} + gV\frac{1}{2}\frac{d}{dt}u^{2}$$
$$= \frac{d}{dt}guu_{t} - g_{t}uu_{t} -gu_{t}^{2} - \frac{d}{dx}(guu_{x}) +gu_{x}^{2} + \frac{1}{2}\frac{d}{dt}(gVu^{2}) - \frac{1}{2}g_{t}Vu^{2}$$
$$= \frac{d}{dt}(guu_{t}) - \frac{1}{2}\frac{d}{dt}(g_{t}u^{2}) + \frac{1}{2}g_{tt}u^{2} - gu_{t}^{2} - \frac{d}{dx}(guu_{x}) + gu_{x}^{2} + \frac{1}{2}\frac{d}{dt}(gVu^{2}) - \frac{1}{2}g_{t}Vu^{2} = 0. $$
{\bf Step 3.} Multiplying the both sides of (1.1) by $hu_{x}$ and rearranging it one obtains:\\
$$hu_{x}u_{tt} - hu_{x}u_{xx} + hu_{x}Vu_{t}$$
$$= h(\frac{d}{dt}u_{x}u_{t}) - h\frac{d}{dt}(u_{x})u_{t} - \frac{1}{2}h\frac{d}{dx}(u_{x}^{2}) + hVu_{t}u_{x}$$
$$= \frac{d}{dt}(hu_{x}u_{t}) - h_{t}u_{x}u_{t} - \frac{1}{2}h\frac{d}{dx}(u_{t}^{2}) - \frac{1}{2}\frac{d}{dx}(hu_{x}^{2}) + \frac{1}{2}h_{x}u_{x}^{2} + hVu_{t}u_{x}$$
$$= \frac{d}{dt}(hu_{t}u_{x}) - h_{t}u_{x}u_{t} - \frac{1}{2}\frac{d}{dx}(hu_{t}^{2}) + \frac{1}{2}h_{x}u_{t}^{2} - \frac{1}{2}\frac{d}{dx}(hu_{x}^{2}) + \frac{1}{2}h_{x}u_{x}^{2} + hVu_{t}u_{x} = 0.$$
{\bf Step 4.} Add all final identities from Step 1 to Step 3 and integrate over $[0,\infty)$. Then, it follows from the boundary condition (1.3) that one can get:
$$\frac{1}{2}\frac{d}{dt}\int^{\infty}_{0}[f(u_{t}^{2} + u_{x}^{2}) + 2guu_{t} + (gV - g_{t})u^{2} + 2hu_{t}u_{x}]dx$$
$$+ \frac{1}{2}\int^{\infty}_{0}(2fV - f_{t} - 2g + h_{x})u_{t}^{2}dx + \frac{1}{2}\int^{\infty}_{0}(2g -f_{t} + h_{x})u_{x}^{2}dx$$
$$+ \frac{1}{2}\int^{\infty}_{0}(g_{tt} - g_{t}V)u^{2}dx +\int^{\infty}_{0}(hV - h_{t})u_{x}u_{t}dx - \frac{1}{2}\int^{\infty}_{0}\frac{d}{dx}(hu_{x}^{2})dx = 0,$$
which implies the desired identity.\ \ $_{\Box}$
\vspace{3mm}

Since the finite speed of propagation property can be applied to the solution $u(t,x)$ of the corresponding problem $(1.1)-(1.3)$, in order to estimate the functions ${\cal E}(t)$ and ${\cal F}(t)$ it suffices to consider all the spatial integrand over the closed interval $[0,R + t]$ ($t \geq 0$).

Now, because of the assumption ({\bf A}) such as $V(L_{2}) > 0$ and the continuity of $V(x)$, there exist constants $L_{1} \in [0,L_{2})$ and $V_{m} > 0$ such that
\[V(x) \geq V_{m}\quad (x \in [{L_{1},L_{2}}]).\]
By using these constants we choose the functions $f(t)$, $g(t)$ and $h(t,x)$ as follow:
$$f(t) = \epsilon_{1}(1+t)^{2},\ \ \ \ g(t) = \epsilon_{2}(1+t),\ \ \ h(t,x) = \epsilon_{3}(1+t)\phi(x),$$
where the {\bf monotone increasing} function $\phi \in C^{\infty}([0,\infty))$ can be defined to satisfy 

\begin{equation}\phi(x) = \left\{
  \begin{array}{cc}
   \displaystyle (1+x) &
        (0 \leq x < L_{1}), \\[3mm]
   \displaystyle 1+L_{2} &
        (L_{2} \leq x),
   \end{array} \right. \end{equation}
where $\epsilon_{i} > 0 \ \ (i=1,2,3)$ are some constants determined later on.
\begin{lem} Let $f$, $g$ and $h$ be defined by $(2.1)$. Then one has the following estimates:\,for each $t \geq 0$,\\
{\rm (i)}\,$\displaystyle{\int^{\infty}_{0}}(hV - h_{t})u_{t}u_{x}dx \geq -\displaystyle{\frac{k}{2}}\displaystyle{\int^{\infty}_{0}}hVu_{t}^{2}dx -\displaystyle{\frac{1}{2k}}\displaystyle{\int^{\infty}_{0}}hVu_{x}^{2}dx -\displaystyle{\frac{1}{2}}\displaystyle{\int^{\infty}_{0}}h_{t}u_{t}^{2}dx -\displaystyle{\frac{1}{2}}\displaystyle{\int^{\infty}_{0}}h_{t}u_{x}^{2}dx,$\\
{\rm (ii)}\,$\displaystyle{\int^{\infty}_{0}\frac{d}{dx}}(hu_{x}^2)dx \leq 0,$\\
where $k > 0$ is a constant.
\end{lem}
{\it Proof.}\,(i): This is easily derived from the following inequality with some positive number $p$:\\
\[\vert u_{t}u_{x}\vert \leq \frac{1}{2p}u_{t}^{2} + \frac{p}{2}u_{x}^{2}.\]
(ii): Because of the finite speed of propagation property, one has
$$\int^{\infty}_{0}\frac{d}{dx}(hu_{x}^{2})dx = h(t,\infty)u_{x}(t,\infty)^{2} - h(t,0)u_{x}(t,0)^{2}$$
$$=-h(t,0)u_{x}(t,0)^{2} \leq 0.\ \ \ _{\Box}$$

Thus, it follows from Lemmas 2.2 and 2.3 that for $t \geq 0$ 
$$\frac{1}{2}\frac{d}{dt}\int^{\infty}_{0}[f(u_{t}^{2} + u_{x}^{2}) + 2guu_{t} + (gV - g_{t})u^{2} + 2hu_{t}u_{x}]dx$$
$$+\frac{1}{2}\int^{\infty}_{0}(2fV - f_{t} - 2g + h_{x} -khV - h_{t})u_{t}^{2}dx$$
\begin{equation}
+ \frac{1}{2}\int^{\infty}_{0}(2g -f_{t} + h_{x} -\frac{1}{k}hV - h_{t})u_{x}^{2}dx + \frac{1}{2}\int^{\infty}_{0}(g_{tt}-g_{t}V)u^{2}dx \leq 0.
\end{equation}
Furthermore, one can get the estimates below:

\begin{lem}
Let f, g and h be defined by $(2.1)$, and assume {\rm (}{\bf A}{\rm )}. If all parameters $\epsilon_{i}$\,{\rm (}$i = 1,2,3${\rm )} are well-chosen, then there exists a large $t_{0} > 0$ such that for all $t \geq t_{0} \gg 0,$ one has\\
{\rm (iii)}\,$2f(t)V(x) - f_{t}(t) -2g(t) + h_{x}(t,x) - kh(t,x)V(x) - h_{t}(t,x) > 0, \quad x \in [0,\infty)$,\\
{\rm (iv)}\,$2g(t) - f_{t}(t) + h_{x}(t,x) - \displaystyle{\frac{1}{k}}h(t,x)V(x) - h_{t}(t,x) > 0, \quad x \in [0,\infty)$\\
with some constant $k > 0$.
\end{lem}
{\it Proof.}\,We first check (iii) by separating the integrand of $x$ into $3$ parts $[0,L_{1}]$, $[L_{1},L_{2}]$, and $[L_{2},+\infty)$.\\
({\rm iii})\,The case $\underline{0 \leq x \leq L_{1}}$:\\
$$2fV - f_{t} - 2g + h_{x} - khV - h_{t}$$
$$= 2\epsilon_{1}(1+t)^{2}V(x) - 2\epsilon_{1}(1+t) - 2\epsilon_{2}(1+t) +\epsilon_{3}(1+t) -k\epsilon_{3}(1+t)(1+x)V(x) - \epsilon_{3}(1+x)$$
$$\geq (1+t)^{2}\{2\epsilon_{1} - k\epsilon_{3}\frac{1+L_{1}}{1+t}\}V(x) + (1+t)\{\epsilon_{3} -2\epsilon_{1} -2\epsilon_{2}  - \epsilon_{3}\frac{1+L_{1}}{1+t}\}.$$

The case $\underline{L_{1} \leq x \leq L_{2}}$:\\
$$2fV - f_{t} - 2g + h_{x} - khV - h_{t}$$
$$= 2\epsilon_{1}(1+t)^{2}V(x) - 2\epsilon_{1}(1+t) - 2\epsilon_{2}(1+t) + \epsilon_{3}(1+t)\phi'(x) - k\epsilon_{3}(1+t)\phi(x)V(x) - \epsilon_{3}\phi(x)$$
$$\geq (1+t)^{2}\{2\epsilon_{1}V_{m} - 2\epsilon_{1}\frac{1}{1+t} - 2\epsilon_{2}\frac{1}{1+t} - k\epsilon_{3}V_{M}\frac{1+L_{2}}{1+t} -\epsilon_{3}\frac{1+L_{2}}{(1+t)^{2}}\},$$

where $V_{M}:=\sup\{V(x)|0\leq x \leq L_{2}\} > 0$, and we have just used the fact that $\phi'(x) \geq 0$.\\

The case $\underline{L_{2} \leq x}$:\, Here we use the finite speed of propagation property of the solution below:\\
$$2fV - f_{t} - 2g + h_{x} - khV - lh_{t}$$
$$= 2\epsilon_{1}(1+t)^{2}V(x) - 2\epsilon_{1}(1+t) - 2\epsilon_{2}(1+t) - k\epsilon_{3}(1+t)(1+L_{2})V(x) - \epsilon_{3}(1+L_{2})$$
$$\geq 2\epsilon_{1}(1+t)^{2}\frac{V_{0}}{1+x} - 2\epsilon_{1}(1+t) - 2\epsilon_{2}(1+t) - k\epsilon_{3}(1+t)(1+L_{2})\frac{V_{1}}{1+x} - \epsilon_{3}(1+L_{2})$$
$$\geq \frac{(1+t)^{2}}{1+x}\{2\epsilon_{1}V_{0} - 2\epsilon_{1}\frac{1+R+t}{1+t} - 2\epsilon_{2}\frac{1+R+t}{1+t} - k\epsilon_{3}(1+L_{2})\frac{V_{1}}{1+t} - \epsilon_{3}(1+L_{2})\frac{1+R+t}{(1+t)^{2}} \}.$$
Note that in this region, $\phi'(x) = 0$.
\vspace{5mm}

Next one can check the condition (iv) similarly to the case (iii).

(iv)\,The case $\underline{0 \leq x \leq L_{1}}$:\\
$$2g - f_{t} + h_{x} - \frac{1}{k}hV - h_{t}$$
$$= 2\epsilon_{2}(1+t) - 2\epsilon_{1}(1+t) + \epsilon_{3}(1+t) - \frac{\epsilon_{3}}{k}(1+t)(1+x)V(x) - \epsilon_{3}(1+x)$$
$$\geq (1+t)\{2\epsilon_{2} - 2\epsilon_{1} + \epsilon_{3} - \frac{\epsilon_{3}}{k}(1+L_{1})V_{M} - \epsilon_{3}\frac{1+L_{1}}{1+t}\}. $$

The case $\underline{L_{1} \leq x \leq L_{2}}$:\\
$$2g - f_{t} + h_{x} - \frac{1}{k}hV - h_{t}$$
$$= 2\epsilon_{2}(1+t) - 2\epsilon_{1}(1+t) + \epsilon_{3}(1+t)\phi'(x) - \frac{\epsilon_{3}}{k}(1+t)\phi(x)V(x) - \epsilon_{3}\phi(x)$$
$$\geq (1+t)\{2\epsilon_{2} - 2\epsilon_{1} - \frac{\epsilon_{3}}{k}(1+L_{2})V_{M} - \epsilon_{3}\frac{1+L_{2}}{1+t}\}.$$

The case $\underline{L_{2} \leq x}$:\\
$$2g - f_{t} + h_{x} - \frac{1}{k}hV - h_{t}$$
$$= 2\epsilon_{2}(1+t) - 2\epsilon_{1}(1+t) - \frac{\epsilon_{3}}{k}(1+t)(1+L_{2})V(x) - \frac{\epsilon_{3}}{l}(1+L_{2})$$
$$\geq (1+t)\{2\epsilon_{2} - 2\epsilon_{1} - \frac{\epsilon_{2}}{k}(1+L_{2})\frac{V_{1}}{1+x} - \epsilon_{2}\frac{1+L_{2}}{1+t}\}$$
$$\geq (1+t)\{2\epsilon_{2} - 2\epsilon_{1} - \frac{\epsilon_{3}}{k}V_{1} - \epsilon_{2}\frac{1+L_{2}}{1+t}\}.$$

So, in order to obtain (iii) and (iv) for large $t \gg 1$ the following five conditions should be satisfied to $\epsilon_{i}$ ($i = 1,2,3$) and $k > 0$:
\begin{equation}
\epsilon_{3} - 2\epsilon_{1} - 2\epsilon_{2} > 0,
\end{equation}

\begin{equation}
2\epsilon_{1}V_{0} - 2\epsilon_{1} - 2\epsilon_{2} > 0,
\end{equation}

\begin{equation}
2\epsilon_{2} - 2\epsilon_{1} + \epsilon_{3} - \frac{\epsilon_{3}}{k}(1+L_{1})V_{M} > 0,
\end{equation}

\begin{equation}
2\epsilon_{2} - 2\epsilon_{1} - \frac{\epsilon_{3}}{k}(1+L_{2})V_{M} > 0,
\end{equation}

\begin{equation}
2\epsilon_{2} - 2\epsilon_{1} - \frac{\epsilon_{3}}{k}V_{1} > 0.
\end{equation}
(2.3) and (2.4) come from the check of (iii), while (2.5)-(2.7) have its origin in the case when we check (iv) as $t \to +\infty$.
  
We need to look for the constants $\epsilon_{1}$, $\epsilon_{2}$, $\epsilon_{3}$, $k > 0$ satisfying these five conditions (2.3)-(2.7) above.\\

First, the condition $(2.5)$ and $(2.6)$ and $(2.7)$ can be unified by the following one condition $(2.8)$:
$$2\epsilon_{2} - 2\epsilon_{1} + \epsilon_{3} -\frac{\epsilon_{3}}{k}(1+L_{1})V_{M} \geq 2\epsilon_{2} - 2\epsilon_{1} - \frac{\epsilon_{3}}{k}(1+L_{2})V_{M}$$
\begin{equation}
\geq 2\epsilon_{2} - 2\epsilon_{1} - \frac{\epsilon_{3}}{k}V^{*} > 0,
\end{equation}
where
\[V^{*} := \max\{(1+L_{2})V_{M},\,V_{1}\} > 0.\]
Thus it suffices to check only $3$ conditions $(2.3)$, $(2.4)$ and $(2.8)$ above. However, for its purpose it is enough to choose all constants $\epsilon_{i} > 0$ ($i = 1,2,3$) and $k > 0$ as follows:
\[\epsilon_{1} := 1,\quad \epsilon_{2} := \frac{V_{0}}{2},\quad\epsilon_{3} := 2V_{0},\]
and
\[k:= \frac{4V_{0}V^{*}}{V_{0}-2},\]
where the assumption $V_{0} > 2$ is essentially used. Therefore, one has the desired estimates if one takes $t \gg 1$ sufficiently large. \ \ \ \ $_{\Box}$ \\

Next lemma is a direct consequence of (2.2) and Lemma 2.4.

\begin{lem}\,Let $u \in X_{1}(0,\infty)$ be the solution to problem $(1.1)-(1.3)$, and $f, g, h$ be defined by $(2.1)$. Under the condition $({\bf A})$, the following estimate holds true:
\begin{equation}\frac{d}{dt}\{f(t)E(t) + g(t)(u_{t},u) + 2(hu_{x},u_{t})\} \leq \frac{1}{2}\frac{d}{dt}\int^{\infty}_{0}(g_{t} - gV)u^{2}dx + \frac{1}{2}\int^{\infty}_{0}(g_{t}V - g_{tt})u^{2}dx\end{equation}
for all $t \geq t_{0}$, where $t_{0} \gg 1$ is a fixed time defined in Lemma {\rm 2.4}.
\end{lem}

Then we integrate the both sides of $(2.9)$ over $[t_{0},t]$:
$$f(t)E(t) + g(t)(u_{t},u) + 2(hu_{x},u_{t})$$
$$\leq f(t_{0})E(t_{0}) + g(t_{0})(u_{t}(t_{0}),u_{t}(t_{0})) + 2(h(t_{0})u_{x}(t_{0}),u_{t}(t_{0}))$$
$$+ \frac{1}{2}\int^{\infty}_{0}(g_{t}-gV)u^{2}dx - \frac{1}{2}\int^{\infty}_{0}(g_{t}(t_{0}) - g(t_{0})V)u^{2}(t_{0})dx + \frac{1}{2}\int^{t}_{t_{0}}\int^{\infty}_{0}(g_{t}V - g_{tt})u^{2}dxds$$
\begin{equation}= C + \frac{1}{2}\int^{\infty}_{0}(g_{t} - gV)u^{2}dx + \frac{1}{2}\int^{t}_{t_{0}}\int^{\infty}_{0}(g_{t}V - g_{tt})u^{2}dxds,\end{equation}
where the constant $C > 0$, which is independent from $t \geq t_{0}$, is defined by
$$C := f(t_{0})E(t_{0}) + g(t_{0})(u_{t}(t_{0}),u_{t}(t_{0})) + 2(h(t_{0})u_{x}(t_{0}),u_{t}(t_{0})) $$
$$- \frac{1}{2}\int^{\infty}_{0}(g_{t}(t_{0}) - g(t_{0})V)u^{2}(t_{0})dx.$$

Furthermore, we have the following estimates.

\begin{lem}\,Let $g$ be defined by $(2.1)$. The the smooth function $g$ satisfies the following two estimates:
$$(v)\ \ g_{t} - gV(x) \leq C_{1},\ \ \ \ x \in [0,\infty),\quad t \geq 0,$$
$$(vi)\ \ g_{t}V - g_{tt} \leq C_{2}V(x),\ \ \ \ x \in [0,\infty),\quad t \geq 0,$$
where $C_{i} > 0\ \ (i = 1,2)$ are some constants.
\end{lem}
{\it Proof.}\ \ \ The proof can be easily checked, so we omit it.\ \ \ $_{\Box}$

On the other hand, we shall prepare the following crucial lemma based on the one dimensional Hardy-Sobolev inequality in the half space case, which is stated in Lemma 2.1.

\begin{lem}\,Let $u \in X_{1}(0,\infty)$ be the solution to problem $(1.1)-(1.3)$. Then it is true that
\begin{equation}
\|u(t,\cdot)\|^{2} + \int^{t}_{0}\int^{\infty}_{0}V(x)|u(s,x)|^{2}dxds \leq C(\|u_{0}\|^{2} + \|(V(\cdot)u_{0} + u_{1})\|^{2}_{1,1/2}),
\end{equation}
provided that $\|(V(\cdot)u_{0} + u_{1})\|_{1,1/2} < +\infty$.
\end{lem}
{\it Proof.}\, The original idea comes from \cite{IM}. We introduce an auxiliary function
$$w(t,x) := \int^{t}_{0}u(s,x)ds.$$
Then $w(t,x)$ satisfies
\begin{equation}
w_{tt} - w_{xx} + V(x)w_{t} = V(x)u_{0} + u_{1},\ \ \ \ (t,x) \in (0,\infty) \times (0,\infty),
\end{equation}
\begin{equation}
w(0,x) = 0,\ \ \  w_{t}(0,x) = u_{0}(x),\ \ \ \  x \in (0,\infty).
\end{equation}
Multiplying $(2.12)$ by $w_{t}$ and integrating over $[0,t]\times [0,\infty)$ we get
$$\frac{1}{2}(\|w_{t}(t,\cdot)\|^{2} + \|w_{x}(t,\cdot)\|^{2}) + \int^{t}_{0}\|\sqrt{V(\cdot)}w_{s}(s,\cdot)\|^{2}ds$$
\begin{equation}= \frac{1}{2}\|u_{0}\|^{2} + \int^{t}_{0}(V(\cdot)u_{0} + u_{1}, w_{s})ds.
\end{equation}

Next step is to use Lemma 2.1 to obtain a series of inequalities below:\\

$$\int^{t}_{0}(V(\cdot)u_{0} + u_{1},w_{s})ds = \int^{t}_{0}\frac{d}{ds}(V(\cdot)u_{0} + u_{1},w)ds$$
$$\leq \int^{\infty}_{0}\sqrt{1+x}|V(x)u_{0} + u_{1}|\frac{|w(t,x)|}{\sqrt{1+x}}dx$$
$$\leq (\sup_{x\in [0,\infty)}\frac{|w(t,x)|}{\sqrt{1+x}})\|V(\cdot)u_{0} + u_{1}\|_{1,1/2}$$
\begin{equation}\leq \frac{1}{4}\|w_{x}\|^{2} + \|V(\cdot)u_{0} + u_{1}\|^{2}_{1,1/2}.
\end{equation}

Combining $(2.14)$ with $(2.15)$ we can derive

$$\frac{1}{2}\|w_{t}(t,\cdot)\|^{2} + \frac{1}{4}\|w_{x}(t,\cdot)\|^{2} + \int^{t}_{0}\int^{\infty}_{0}V(x)w_{t}(s,x)dxds$$
$$\leq \frac{1}{2}\|u_{0}\|^{2} + \|V(\cdot)u_{0} + u_{1}\|_{1,1/2}^{2}.$$
The desired estimate follows from the estimate above and the fact that $w_{t}=u$.\ \ \ \ \ $_{\Box}$
\vspace{3mm}

It follows from (2.10) and Lemmas 2.6 and 2.7 that there exists a constant $C > 0$ such that
\begin{equation}
f(t)E(t) + g(t)(u(t,\cdot),u_{t}(t,\cdot)) + 2(hu_{x},u_{t}) \leq C \quad (t \geq t_{0}),
\end{equation}
provided that $\|(V(\cdot)u_{0} + u_{1})\|_{1,1/2} < +\infty$. 
 
Finally, we can derive the following lemma.

\begin{lem}\,Let $h$ be defined by $(2.1)$. Then, for all $t \geq t_{0} \gg 1$ it is true that
$$f(t)E(t) + 2(hu_{x},u_{t}) \geq Cf(t)E(t),$$
where $C>0$ is a constant.
\end{lem}
{\it Proof.}\, Indeed, one has
$$f(t)E(t) + 2(hu_{x},u_{t}) \geq \frac{1}{2}\int^{\infty}_{0}f(t)(u_{t}^{2} + u_{x}^{2})dx - \int^{\infty}_{0}h(t,x)(u_{x}^{2} + u_{t}^{2})dx$$
$$\geq \frac{1}{2}\int^{\infty}_{0}(f(t) - 2h(t,x))(u_{x}^{2} + u_{t}^{2})dx.$$

On the other hand, if necessarily, by choosing $t_{0} \gg 1$ further large enough, one can derive the following estimates for $t \geq t_{0}$:
$$f(t) - 2h(t,x) \geq \epsilon_{1}(1+t)^{2} - 2\epsilon_{3}(1+t)(1+L_{2})$$
$$= (1+t)^{2}\{\epsilon_{1} - \frac{2\epsilon_{3}(1+L_{2})}{1+t}\} \geq C(1+t)^{2} \geq Cf(t),$$
with some constant $C > 0$. Here we have just used the monotonicity of the function $\phi(x)$ closely related with the definition of the function $h(t,x)$.\ \ \ \ \ $_{\Box}$

Now we can finalize the proof of Theorem 1.1.\\
{\it Proof of Theorem 1.1.}\ \ \ We first note that one can use Lemma 2.7 because one can check $\|(V(\cdot)u_{0} + u_{1})\|_{1,1/2} < +\infty$ under the assumption on the initial data stated in Theorem 1.1. Thus, 
by using the Schwarz inequality, (2.16) and Lemma 2.8 we get

$$C(1+t)^{2}E(t) \leq g(t)\|u(t,\cdot)\|\|u_{t}(t,\cdot)\| + C \leq Cg(t)\sqrt{E(t)} + C,\ \ \ t \geq t_{0}.$$

Furthermore, if we set $X(t)=\sqrt{E(t)}$ for $t \in [0,+\infty)$, then one has

\begin{equation}f(t)X(t)^{2} - Cg(t)X(t) - C \leq 0,\ \ \ \ t \geq t_{0}. \end{equation}

By solving the quadratic inequality $(2.17)$ for $X(t)$ we have

$$\sqrt{E(t)} \leq \frac{Cg(t) + \sqrt{C^{2}g(t)^{2} + 4Cf(t)}}{2f(t)}\quad (t \geq t_{0}).$$

This inequality leads to
$$E(t) \leq C\biggl(\frac{g(t)}{f(t)}\biggr)^{2} + C\biggl(\frac{1}{f(t)}\biggr),\ \ \ t \geq t_{0},$$

which implies the desired decay estimates.\ \ \ \ $_{\Box}$



\par
\vspace{0.2cm}
\par
\noindent{\em Acknowledgment.}
\smallskip
The work of the first author (R. IKEHATA) was supported in part by Grant-in-Aid for Scientific Research (C) 15K04958 of JSPS.



\begin{thebibliography}{99}
\bibitem{AIK} L. Aloui, S. Ibrahim and M. Khenissi, Energy decay for linear dissipative wave equations in exterior domains, J. Diff. Eqns 259 (2015), 2061-2079.
\bibitem{D} M. Daoulatli, Energy decay rates for solutions of the wave equation with linear damping in exterior domain, arXiv:1203.6780v4.
\bibitem{ike-1} R. Ikehata, A remark on a critical exponent for the semilinear dissipative wave equation in the one dimensional half space, Diff. Int. Eqns 16, No. 6 (2003), 727-736.
\bibitem{ike-2} R. Ikehata, Fast decay of solutions for linear wave equations with dissipation localized near infinity in an exterior domain, J. Diff. Eqns 188 (2003), 390-405.
\bibitem{ike-3} R. Ikehata and Y. Inoue, Total energy decay for semilinear wave equations with a critical potential type of damping, Nonlinear Anal. 69 (2008), 1396-1401.
\bibitem{IM} R. Ikehata and T. Matsuyama, $L^{2}$-behaviour of solutions to the linear heat and wave equations in exterior domains, Sci. Math. Japon. 55 (2002), 33-42.
\bibitem{ike-4} R. Ikehata, G. Todorova and B. Yordanov, Optimal decay rate of the energy for linear wave equations with a critical potential, J. Math. Soc. Japan 65, No. 1 (2013), 183-236.

\bibitem{L} G. Lebeau, \'Equations des ondes amorties, Algebraic and geometric methods in Math. Physics, A. Boutet de Monvel and V. Marchenko (eds), Kluwer Academic, The Netherlands, (1996), 73-109.
\bibitem{Ma} A. Matsumura, On the asymptotic behavior of solutions of semilinear wave equations, Publ. RIMS Kyoto Univ. 12 (1976), 169-189.

\bibitem{Mochi} K. Mochizuki, Scattering theory for wave equations with dissipative terms, Publ. Res. Inst. Math. Sci. 12 (1976), 383-390.
\bibitem{MN} K. Mochizuki and H. Nakazawa, Energy decay and asymptotic behavior of solutions to the wave equations with linear dissipation, Publ. Res. Inst. Math. Sci. 32 (1996), 401-414.
\bibitem{Mora} C. Morawetz, The decay of solutions of the exterior initial-boundary value problem for the wave equation, Comm. Pure Appl. Math. 14 (1961), 561-568.

\bibitem{N} M. Nakao, Energy decay for the linear and semilinear wave equations in exterior domains with some localized dissipations, Math. Z. 238 (2001), 781-797.


\bibitem{Ni} K. Nishihara, $L^{p}$-$L^{q}$ estimates to the damped wave equation in $3$-dimensional space and their application, Math. Z. 244 (2003), 631-649.
\bibitem{R} R. Racke, Decay rates for solutions of damped systems and generalized Fourier transforms, J. Reine Angew. Math. 412 (1990), 1-19.
\bibitem{TY} G. Todorova and B. Yordanov, Weighted $L^{2}$-estimates of dissipative wave equations with variable coefficients, J. Diff. Eqns 246 (2009), 4497-4518.

\bibitem{U} H. Uesaka, The total energy decay of solutions for the wave equation with a dissipative term, J. Math. Kyoto Univ. 20 (1980), 57-65.
\bibitem{W} T. Watanabe, Global existence and decay estimates for quasilinear wave equations with nonuniform dissipative term, Funk. Ekvac. 58 (2015), 1-42.

\bibitem{Z} E. Zuazua, Exponential decay for the semilinear wave equation with localized damping in unbounded domains, J. Math. Pures Appl. 70 (1991), 513-529.

\end{thebibliography}
\end{document}